\documentclass[12pt]{article}
\usepackage[english]{babel}
\usepackage{amsmath,amsfonts,amssymb,amsthm} 

\usepackage{epsf}
\usepackage{graphicx} 

\usepackage[numbers]{natbib}

\usepackage{float}
\usepackage{tikz}
\usetikzlibrary{positioning, calc, fit, decorations.pathreplacing}

\usepackage{hyperref} 
\usepackage{mathrsfs} 

\newtheorem{claim}{Claim}

\newcommand\restr[2]{{
  \left.\kern-\nulldelimiterspace 
  #1 
  \vphantom{\big|} 
  \right|_{#2} 
  }} 

\newcommand\subsetsim{\mathrel{\substack{
  \textstyle\subset\\[-0.2ex]\textstyle\sim}}}
  
\renewcommand{\vec}[1]{\mathbf{#1}}
\newcommand{\mtrxN}[1]{#1}
\newcommand{\mtrxNN}[1]{\mathsf{#1}}
\newcommand{\mtrxNNs}[1]{\mathcal{#1}}
\newcommand{\operator}[1]{\mathscr{#1}}

\title{Spectral analysis of a class of IRK stage preconditioners for linear wave equations}
\author{Michal Outrata\footnote{Charles University, Czechia, e-mail: outrata@karlin.mff.cuni.cz \\ This work was supported by the grants UNCE/24/SCI/005 and PRIMUS/25/SCI/022 of Charles University in Prague.} and Aman Rani\footnote{Johns Hopkins University, MD, USA, e-mail: arani1@jh.edu}}
\date{}

\begin{document}
\maketitle

\begin{abstract}
We lay out the analysis of the ``Kronecker preconditioners'' adopted in~\cite{rani2025efficient} for the linear wave equations. Importantly, we deliver the analysis justifying the claims in~\cite[Section 4.2]{rani2025efficient} about these preconditioners and we unpack all of the necessary numerical background for a full understanding of the observed phenomena.
\end{abstract}
 

\paragraph*{MSC2020 class}
65F08, 65L06, 65M22

\section{Introduction and background}\label{sec_IntroBckgrnd}
In~\cite{rani2025efficient}, the authors proposed a new approach for accelerating the solution process for the class of linear wave equations of the form
\begin{equation}\label{eqn_secIntro_PDE_utt_m_L_eq_f}
\left( \partial_{tt} - \operator{L} \right) u = f \quad \mathrm{in} \; \Omega \times (0,T)
\end{equation}	
\noindent with some initial and boundary conditions $u^{(0)}(\vec{x})$ and $g(\vec{x},t)$, where $(0,T)$ and $\Omega$ are temporal and spatial domains of interest, $\operator{L}$ is a linear, second-order spatial operator and $f = f(\vec{x},t)$ is a source function. While the canonical model problem of \emph{the} wave equation corresponds to $\operator{L}=c^2\Delta$ with the wave speed $c^2>0$, the analysis applies to a much wider group of $\operator{L}$, see Section~\ref{sec_precGMRESbeh}. The approach relies on a reformulation of~\eqref{eqn_secIntro_PDE_utt_m_L_eq_f} into a larger, first-order system for which the authors use standard finite-elements discretization in space and (fully) implicit Runge-Kutta methods in time, proposing new, efficient preconditioners for the resulting linear systems. Generally speaking, using (fully) implicit Runge-Kutta methods in this context and especially proposing and analyzing suitable preconditioners for the resulting system of algebraic equations has received a renewed interest in the last couple of years, see~\cite{howle2021new,farrell2021irksome,howle2022efficient,southworth2022fastI,southworth2022fastII,dravins2024stageparallel_2,
dravins2024performance} and also~\cite{outrata2022IRKILAS,outrata2023irksstage,dravins2024stageparallel_1,dravins2024spectral,pearson2024pint,outrata2025recent}. 

The purpose of this work is to provide theoretical analysis for the performance of the preconditioners considered in~\cite{rani2025efficient} based on the Sylvester reformulation, generalizing the analysis in~\cite{outrata2023irksstage}. We refer the interested reader to the aforementioned works for more exhaustive literature overview. We first summarize the background and our notation in Section~\ref{sec_IntroBckgrnd}, then derive the spectral analysis in Section~\ref{sec_SpctrlAnal} and finally make use of it for estimating the GMRES convergence for model problems in Section~\ref{sec_precGMRESbeh}.

\subsection{Discretization and implicit Runge-Kutta methods}
For the sake of space we omit several discretization details, following the set-up in~\cite[Section 2]{rani2025efficient}. We introduce an auxiliary function $w(\vec{x},t) = [u(\vec{x},t), u_t(\vec{x},t)]^T$ and rewrite~\eqref{eqn_secIntro_PDE_utt_m_L_eq_f} as
\begin{equation}\label{eqn_secIntro_PDE_wt_eq_B_p_f}
w_t = 
\begin{bmatrix} & \mathscr{I} \\ \mathscr{L}
\end{bmatrix} w + \begin{bmatrix} \\ f \end{bmatrix}  \quad \mathrm{in} \; (0,T)\times \Omega,
\end{equation}
\noindent equipped with appropriate initial and boundary conditions. We discretize $\mathscr{L}$ using a simple FEM scheme with linear Lagrangian polynomials on a conforming triangular mesh (but we do not build on any properties specific to P1-FEM or FEM in general), obtaining the mass and stiffness matrices $M \in \mathbb{R}^{N\times N}$ and $E \in \mathbb{R}^{N\times N}$ so that instead of~\eqref{eqn_secIntro_PDE_wt_eq_B_p_f} we solve the approximate system of first order ordinary differential equations
\begin{equation}\label{eqn_secIntro_PDE_DiscrtzdInSpc_SysOfFrstOrdrODEs}
\mtrxNN{M} \vec{w}_t = \mtrxNN{B} \vec{w} + \vec{b}  \quad \mathrm{in} \; (0,T),
\end{equation}
\noindent where
\begin{equation*}
\mtrxNN{M} = \begin{bmatrix} \mtrxN{M} \\ & \mtrxN{M} \end{bmatrix}
\quad \mathrm{and} \quad
\mtrxNN{B} = \begin{bmatrix} & \mtrxN{M} \\ -\mtrxN{E} \end{bmatrix},
\end{equation*}
\noindent and $\vec{b} \in \mathbb{R}^{2N}$ aggregates the appropriate terms coming from the discretization of $f$ and the other conditions imposed on~\eqref{eqn_secIntro_PDE_utt_m_L_eq_f} (including boundary conditions) and we impose the correct initial condition based on~\eqref{eqn_secIntro_PDE_wt_eq_B_p_f}.
We solve~\eqref{eqn_secIntro_PDE_DiscrtzdInSpc_SysOfFrstOrdrODEs} numerically using an $s$-stage implicit Runge-Kutta method. That is, we compute a sequence of approximations $\vec{w}^{(1)} ,\vec{w}^{(2)}, \dotsc$ at timepoints $t_0=0 < t_1=h_t < t_2 = 2h_t < \dotsc < T$, given by
\begin{equation*}
\mtrxNN{M}\vec{w}^{(m)} = \mtrxNN{B}\vec{w}^{(m-1)} + h_t \sum\limits_{i=1}^{s} b_i\vec{k}_i^{(m)},
\end{equation*}
\noindent where $b_1,\dotsc, b_s \in \mathbb{R}$ are the so-called Butcher weights and the vectors $\vec{k}_1^{(m)}, \dotsc , \vec{k}_s^{(m)}$ are the so-called stage vectors. Computing these is \emph{the} bottleneck of the calculation as these are governed by the system
\begin{equation}\label{eqn_secIntro_IRK_StageEqn}
\mtrxNNs{T} \; \vec{k}^{(m)} = \vec{b}^{(\mathrm{IRK})},
\quad \mathrm{with} \quad \mtrxNNs{T} := \left( I_s \otimes \mtrxNN{M} - h_t A \otimes \mtrxNN{B} \right),
\end{equation}
\noindent where $\vec{k}^{(m)} = [(\vec{k}^{(m)}_1)^T,\dotsc,(\vec{k}^{(m)}_s)^T]^T$, $A \in \mathbb{R}^{s\times s}$ is the Butcher matrix, $\otimes$ stands for the Kronecker product and $\vec{b}^{(\mathrm{IRK})}$ aggregates several terms, see~\cite[equation (2.5)]{rani2025efficient}.

\subsection{Solving the linear system}
The problem~\eqref{eqn_secIntro_IRK_StageEqn} is solved almost exclusively by Krylov subspace methods, see~\cite{liesen2013krylov}, and it seems that the GMRES method is the most common choice (although, e.g., in~\cite{Neytcheva2020} the authors use GCR). Either way, a good preconditioner is essential. The standard approach has been to replace the dense matrix $A$ with $\tilde{A} = (\alpha_{ij})_{i,j=1:s}$ so that the solve becomes cheaper -- classical choices include $\tilde{A}$ to be diagonal or triangular. For~\eqref{eqn_secIntro_IRK_StageEqn} we obtain the preconditioner
\begin{equation}\label{eqn_secSolv_Preconditioner}
\mtrxNNs{P}:= I_s \otimes \mtrxNN{M} - h_t \tilde{A} \otimes \mtrxNN{B} ,
\end{equation}
\noindent see also~\cite{staff2006preconditioning} in addition to the cited literature. This manuscript deals with the analysis of the observed performance of such preconditioners in~\cite{rani2025efficient}; there the authors also introduced a new choice for $\tilde{A}$, leading to the so-called \emph{triangular approximate inverse preconditioner} (TAI; see~\cite[Section 4]{rani2025efficient})\footnote{Even though the TAI preconditioner was a new contribution of~\cite{rani2025efficient}, the overall role of the preconditioners of the form~\eqref{eqn_secSolv_Preconditioner} in the manuscript is a ``industry benchmark'', against which the main novelty -- the so-called Sylvester reformulation of the large system -- is measured.}.

We note that an increasingly common practice in the field is to transform the system~\eqref{eqn_secIntro_IRK_StageEqn} by factoring out the mass matrix on the left and the Butcher matrix on the right, obtaining 
\begin{equation*}
\left( I_s \otimes \mtrxNN{M} \right) \left( A^{-1} \otimes \mtrxNN{I}_{2N} - h_t I_s \otimes \mtrxNN{M}^{-1} \mtrxNN{B} \right) \left( A \otimes I_{2N} \right)  \vec{k}^{(m)} = \vec{b}^{(\mathrm{IRK})},
\end{equation*}
\noindent so that the reformulated system looks like
\begin{equation}\label{eqn_secSolv_PreconditionerReformulated}
\left( A^{-1} \otimes \mtrxNN{I}_{2N} - h_t I_s \otimes \mtrxNN{M}^{-1} \mtrxNN{B} \right)  \hat{\vec{k}}^{(m)} = \hat{\vec{b}}^{(\mathrm{IRK})},
\end{equation}
\noindent with $\hat{\vec{b}}^{(\mathrm{IRK})} = \left( I_s \otimes \mtrxNN{M}^{-1} \right) \vec{b}^{(\mathrm{IRK})}$ and $\hat{\vec{k}}^{(m)} = \left( A \otimes \mtrxNN{I}_{2N} \right) \vec{k}^{(m)}$. This idea was first proposed by Butcher in~\cite{butcher1976implementation}; see~\cite[Section 4]{dravins2024stageparallel_1} or~\cite[Section 4.1]{pazner2017stage} for a further discussion. In fact,~\cite[equation (3.1)]{rani2025efficient} corresponds exactly to~\eqref{eqn_secSolv_PreconditionerReformulated}, before the Sylvester reformulation is applied. While the quantitative results of our analysis slightly change depending on whether we consider~\eqref{eqn_secSolv_Preconditioner} or~\eqref{eqn_secSolv_PreconditionerReformulated}, we want to emphasize that the analysis itself is equally applicable to both set-ups and we choose to continue with~\eqref{eqn_secSolv_Preconditioner} only for compatibility reasons with~\cite{rani2025efficient}.

Using GMRES for a (preconditioned) linear system $C\vec{x} = \vec{b}$ and assuming $C$ is diagonalizable, i.e., $C=S\Lambda S^{-1}$ and $\Lambda=\mathrm{diag}(\lambda_1, \dotsc ,\lambda_d)$, the standard \emph{ideal GMRES convergence bound} for the residuals $\vec{r}_{\ell}$ reads
\begin{equation}\label{eqn_secSolv_GMRESPolyInEigenvals}
\frac{\| \vec{r}_{\ell}\|}{\| \vec{r}_0\|} \leq  \; \kappa (S) \; \min\limits_{
\substack{\varphi(0)=1 \\ \mathrm{deg}(\varphi)\leq \ell} }  \max\limits_{1\leq i \leq d} |\varphi(\lambda_i)|,
\end{equation}
\noindent where $\kappa(S)$ is the 2-norm condition number of the matrix $S$, see, e.g.,~\cite[Section 5.7.2]{liesen2013krylov}. As this bound is the starting point of our analysis, we want to first highlight that it can overestimate the GMRES convergence behavior to arbitrary level based on the interaction of $S$, $\vec{x}_0$ and $\vec{b}$. While the spectral information of the system matrix (in our case of the preconditioned one) does not generally govern the convergence (see~\cite{Greenbaum1994a,Greenbaum1996,Arioli1998} and also~\cite[Chapter 2 and Section 5.7]{liesen2013krylov} and references therein) the bound~\eqref{eqn_secSolv_GMRESPolyInEigenvals} separates the contribution of the eigenbasis and spectrum, using the latter for the convergence profile and the former for an initial offset of this profile.

In cases where~\eqref{eqn_secSolv_GMRESPolyInEigenvals} is justifiable, both $\kappa(S)$ and the mixed min-max problem usually need to be further estimated and descriptive estimates of the latter are \emph{difficult} to obtain. Hence, it is common to try simplifying things by replacing the discrete set over which we take the maximum, here $\sigma^{\mathrm{discr}} = \{ \lambda_1, \dotsc , \lambda_d  \}$, by a non-discrete one, which we denote by $\sigma^{\mathrm{non-discr}}$. In general, we aim for $\sigma^{\mathrm{non-discr}}$ such that
\emph{(i)} $\sigma^{\mathrm{discr}} \subsetsim \sigma^{\mathrm{non-discr}}$ and $\sigma^{\mathrm{non-discr}}$ is a finite union of connected, compact sets,
\emph{(ii)} $\partial \sigma^{\mathrm{non-discr}}$ is ``well-captured'' by\footnote{The boundary $\partial \sigma^{\mathrm{non-discr}}$ is to be understood in $\mathbb{C}$ and shows up because the maximum of any polynomial over a connected, bounded set is attained at the boundary of that set, see~\cite[Section 2]{driscoll1998potential}.} $\sigma^{\mathrm{discr}}$, and 
\emph{(iii)} we are able to bound or estimate the simplified problem
\begin{equation}\label{eqn_secSolv_GMRESPolyOverSigmaNonDiscr}
\min\limits_{
\substack{\varphi(0)=1 \\ \mathrm{deg}(\varphi)\leq \ell} }  \max\limits_{ \lambda \in \partial \sigma^{\mathrm{non-discr}} } |\varphi(\lambda)|,
\quad \ell = 1,2,\dotsc .
\end{equation}
\noindent If \emph{(i-ii)} holds, replacing~\eqref{eqn_secSolv_GMRESPolyInEigenvals} with~\eqref{eqn_secSolv_GMRESPolyOverSigmaNonDiscr} introduces only a ``small error'' and vice versa. If $\sigma^{\mathrm{non-discr}}$ is compact and without isolated points, the potential theory allows us to compute the \emph{asymptotic convergence rate}
\begin{equation*}
\rho_{\mathrm{est}} := \lim\limits_{\ell \rightarrow +\infty} \left( \min\limits_{ \substack{\varphi(0)=1 \\ \mathrm{deg}(\varphi)\leq \ell}  } \max\limits_{z \in \sigma^{\mathrm{non-discr}}} |\varphi(z)| \right)^{1/ \ell},
\end{equation*}
\noindent see~\cite{driscoll1998potential,embree1999green}. Having an estimate $\kappa_{\mathrm{est}} \approx \kappa(S)$, we then estimate~\eqref{eqn_secSolv_GMRESPolyInEigenvals} by
\begin{equation}\label{eqn_secSolv_GMRESAsymptConvEst}
\frac{\| \vec{r}_{\ell}\|}{\| \vec{r}_0\|} \lesssim  \; \kappa_{\mathrm{est}} \cdot \rho_{\mathrm{est}}^{\ell},
\quad \ell = 1,2,\dotsc
\end{equation}
\noindent a \emph{linear} convergence \emph{estimate}. Notice that if \emph{(i-ii)} are satisfied reasonably well, then indeed also the min-max part of the bound~\eqref{eqn_secSolv_GMRESPolyInEigenvals} decreases approximately linearly at the rate $\rho_{\mathrm{est}}$. In other words, the failure to capture a potential non-linear GMRES behavior (e.g., an initial phase) is due to the ideal GMRES bound we started with and can be generally attributed to $\kappa(S)$. Altogether, choosing $\sigma^{\mathrm{non-discr}}$ aptly based on \emph{(i--iii)} is key to obtain descriptive estimates of~\eqref{eqn_secSolv_GMRESPolyInEigenvals}; see~\cite[Sections 2 and 3]{outrata2023irksstage} for further details.

\section{Spectral analysis}\label{sec_SpctrlAnal}
Throughout the text we assume that the matrix pencil $\{\mtrxN{M},\mtrxN{E}\}$ is \emph{diagonalizable}, i.e., we assume there exist $N$ linearly independent left and right eigenvectors  $\vec{u}_1,\dotsc ,\vec{u}_N$ and $\vec{v}_1,\dotsc ,\vec{v}_N$ corresponding to the generalized eigenvalues $\lambda_1^{(\mtrxN{M})},\lambda_1^{(\mtrxN{E})},\dotsc ,\lambda_N^{(\mtrxN{M})},\lambda_N^{(\mtrxN{E})}$ so that with
\begin{equation*}
\mtrxN{U} := \left[\vec{u}_1,\cdots ,\vec{u}_N\right] 
\quad \mathrm{and} \quad 
\mtrxN{V} := \left[\vec{v}_1,\dotsc ,\vec{v}_N\right]
\end{equation*}
\noindent and $\Lambda_{\mtrxN{M}} := \mathrm{diag} \left( \lambda_1^{(\mtrxN{M})}, \dotsc, \lambda_N^{(\mtrxN{M})} \right), \Lambda_{\mtrxN{E}} := \mathrm{diag} \left( \lambda_1^{(\mtrxN{E})}, \dotsc, \lambda_N^{(\mtrxN{E})} \right)$, we obtain
\begin{equation*}
\begin{gathered}
\mtrxN{U}^{T} \mtrxN{M} \mtrxN{V} = \Lambda_{\mtrxN{M}}
\quad \mathrm{and} \quad
\mtrxN{U}^{T} \mtrxN{E} \mtrxN{V} = \Lambda_{\mtrxN{E}}.
\end{gathered}
\end{equation*}
\noindent Many important applications lead to such $\mtrxN{M}$ and $\mtrxN{E}$, e.g., if the spatial differential operator $\operator{L}$ is coercive and self-adjoint, then many discretization schemes result in a \emph{symmetric, positive-definite} pencil $\{\mtrxN{M},\mtrxN{E}\}$, i.e., $\mtrxN{E}$ is symmetric and $\mtrxN{M}$ is symmetric, positive-definite. This implies that $\mtrxN{U},\mtrxN{V}$ can be chosen so that $\mtrxN{V}=\mtrxN{U}$ and $\mtrxN{V}^T\mtrxN{M}\mtrxN{V} = \mtrxN{I}$, see~\cite[Sections 2.3, 2.6, 5 and 8]{bai2000templates} for further details.

Notice that both the system matrix~\eqref{eqn_secIntro_IRK_StageEqn} and the preconditioner~\eqref{eqn_secSolv_Preconditioner} are blockwise a linear combination of $\mtrxNN{M}$ and $\mtrxNN{B}$, i.e., both are of the form
\begin{equation}\label{eqn_secAnal_Z_eq_bM_p_cE}
\begin{aligned}
\mtrxNNs{Z} &= 
\begin{bmatrix}
\beta_{11}\mtrxNN{M} + \gamma_{11}\mtrxNN{B} & \cdots & \beta_{1s}\mtrxNN{M} + \gamma_{1s}\mtrxNN{B} \\
\vdots & \ddots & \vdots \\
\beta_{s1}\mtrxNN{M} + \gamma_{s1}\mtrxNN{B} & \cdots & \beta_{ss}\mtrxNN{M} + \gamma_{ss}\mtrxNN{B}
\end{bmatrix} \\
&= 
\begin{bmatrix} \mtrxNN{M} \\ & \ddots \\ & & \mtrxNN{M} \end{bmatrix}
\begin{bmatrix}
\beta_{11}\mtrxNN{I} + \gamma_{11}\mtrxNN{M}^{-1}\mtrxNN{B} & \cdots & \beta_{1s}\mtrxNN{I} + \gamma_{1s}\mtrxNN{M}^{-1}\mtrxNN{B} \\
\vdots & \ddots & \vdots \\
\beta_{s1}\mtrxNN{I} + \gamma_{s1}\mtrxNN{M}^{-1}\mtrxNN{B} & \cdots & \beta_{ss}\mtrxNN{I} + \gamma_{ss}\mtrxNN{M}^{-1}\mtrxNN{B}
\end{bmatrix}.
\end{aligned}
\end{equation}
\noindent Moreover, $\mtrxNN{M}, \mtrxNN{E}$ are both blockwise a linear combination of $\mtrxN{M}$ and $\mtrxN{E}$. Hence, if we can diagonalize the pencil $\{\mtrxN{M},\mtrxN{E}\}$, then this calculation translates to $\mtrxNNs{Z}$. Setting
\begin{equation*}
\mtrxNN{U} :=  
\begin{bmatrix}
\mtrxN{U} \\ & \mtrxN{U}
\end{bmatrix}, \; \,
\mtrxNNs{U} :=  I_s \otimes \mtrxNN{U}
\quad \mathrm{and} \quad 
\mtrxNN{V} :=  
\begin{bmatrix}
\mtrxN{V} \\ & \mtrxN{V}
\end{bmatrix}, \; \,
\mtrxNNs{V} :=  I_s \otimes \mtrxNN{V},
\end{equation*}
\noindent we transform $\mtrxNNs{Z}$ into a block matrix where each block is of the size $2N$-by-$2N$ and is a linear combination of a diagonal (stemming from $\mtrxNN{M}$ terms) and antidiagonal (stemming from $\mtrxNN{B}$ terms) matrices:
\begin{equation*}
\mtrxNNs{U}^T \mtrxNNs{Z} \mtrxNNs{V} = 
\begin{bmatrix} \Lambda_{\mtrxNN{M}} \\ & \ddots \\ & & \Lambda_{\mtrxNN{M}} \end{bmatrix} 
\begin{bmatrix}
\beta_{11} \mtrxNN{I} + \gamma_{11}\Lambda_{\mtrxNN{M}^{-1} \mtrxNN{B}} & \cdots & \beta_{1s}\mtrxNN{I} + \gamma_{1s}\Lambda_{\mtrxNN{M}^{-1} \mtrxNN{B}} \\
\vdots & \ddots & \vdots \\
\beta_{s1}\mtrxNN{I} + \gamma_{s1}\Lambda_{\mtrxNN{M}^{-1} \mtrxNN{B}} & \cdots & \beta_{ss}\mtrxNN{I} + \gamma_{ss}\Lambda_{\mtrxNN{M}^{-1} \mtrxNN{B}}
\end{bmatrix},
\end{equation*}
\noindent where, setting\footnote{Often, $\lambda_k$ is called the $k$-th generalized eigenvalue of the pencil $\{M,E\}$.} $\lambda_k:=\lambda_k^{(E)}/\lambda_k^{(M)}$ and $\Lambda := \mathrm{diag}(\lambda_1,\dotsc , \lambda_N)$, we have
\begin{equation}\label{eqn_secSpctrAnal_Lambda_mtrxNN_def}
\Lambda_{\mtrxNN{M}} = \begin{bmatrix} \Lambda_{\mtrxN{M}} \\ & \Lambda_{\mtrxN{M}} \end{bmatrix}
\quad \mathrm{and} \quad
\Lambda_{\mtrxNN{M}^{-1} \mtrxNN{B}} = \begin{bmatrix} & \mtrxN{I}_{N} \\ -\Lambda_{\mtrxN{M}}^{-1} \Lambda_{\mtrxN{E}} \end{bmatrix} \equiv 
\begin{bmatrix} & \mtrxN{I}_N \\ -\Lambda \end{bmatrix}.
\end{equation}
\noindent This is a very particular structure and can be further simplified by grouping the rows and columns featuring the mode $\lambda_k^{(\mtrxN{M},\mtrxN{E})}$ together for each $k$. For the case $s=2$ this corresponds a well-known \emph{red-black} symmetrical reordering of the equations and the unknowns according. More generally, we set $\pi = [1,N+1,\dotsc, (s-1)N+1 \, , \, 2,N+2,\dotsc, (s-1)N+2 \, ,\dotsc]$ and denote $\Pi \in \mathbb{R}^{2Ns\times 2Ns}$ as the permutation matrix that reorders the columns according to $\pi$. Then, applying the above machinery, we can write
\begin{equation}\label{eqn_secSpctrAnal_BlckDiagonalizationOfmtrxNNsZ_to2s}
\left( \mtrxNNs{U} \Pi \right)^T \mtrxNNs{Z} \mtrxNNs{V} \Pi = 
\begin{bmatrix} \Lambda_{1}^{(\mtrxNN{M})} \\ & \ddots \\ & & \Lambda_{N}^{(\mtrxNN{M})} \end{bmatrix}
\begin{bmatrix} \tilde{Z}_{\lambda_1} \\ & \ddots \\ & & \tilde{Z}_{\lambda_N} \end{bmatrix},
\end{equation}
\noindent and, recalling~\eqref{eqn_secSpctrAnal_Lambda_mtrxNN_def}, we get $\tilde{Z}_{\lambda_k} \in \mathbb{R}^{2s\times 2s}$ as
\begin{equation*}
\tilde{Z}_{\lambda_k} = 
\begin{bmatrix}
\beta_{11} I_{2} + \gamma_{11}\Lambda_{k} & \cdots & \beta_{1s}I_{2} + \gamma_{1s}\Lambda_{k} \\
\vdots & \ddots & \vdots \\
\beta_{s1}I_{2} + \gamma_{s1}\Lambda_{k} & \cdots & \beta_{ss}I_{2} + \gamma_{ss}\Lambda_{k}
\end{bmatrix},
\quad
\begin{aligned}
\Lambda^{(\mtrxNN{M})}_{k} &:= \lambda_k^{(M)} I_2, \\
\Lambda_{k} &:= \begin{bmatrix} & 1 \\ -\lambda_k \end{bmatrix}.
\end{aligned}
\end{equation*}
\noindent Calculating further, we notice that for any $\beta,\gamma \in \mathbb{R}$ and $\lambda\neq 0$ we have
\begin{equation}\label{eqn_secSpctrAnalKroneckr_2by2SclrBlck_ExplDiagonaliztn}
\begin{bmatrix}
\beta  & \gamma  \\
-\gamma \lambda & \beta
\end{bmatrix} = 
F_{\lambda}^{-1}
\begin{bmatrix}
\beta + \vec{i} \gamma \sqrt{\lambda} \\
& \beta - \vec{i} \gamma \sqrt{\lambda}
\end{bmatrix}
F_{\lambda}
\quad \mathrm{with} \quad
F_{\lambda} = \begin{bmatrix}
\vec{i} & \lambda^{-1/2} \\
-\vec{i} & \lambda^{-1/2}
\end{bmatrix}
\end{equation}
\noindent where $F_{\lambda}$ depends only on $\lambda$; in particular, $F_{\lambda}$ is independent of $\beta,\gamma$ and hence~\eqref{eqn_secSpctrAnal_BlckDiagonalizationOfmtrxNNsZ_to2s} can be made more explicit. Considering $\lambda= \lambda_k$ and setting $\tilde{F}_{\lambda_k} := I_s \otimes F_{\lambda_k}$, we observe that $\tilde{F}_{\lambda_k}$ diagonalizes each of the blocks in $\tilde{Z}_{\lambda_k}$ (provided $\lambda_k\neq 0$), i.e.,
\begin{equation*}
\resizebox{.995\textwidth}{!}{$
\tilde{F}_{\lambda_k}^{-1} \tilde{Z}_{\lambda_k} \tilde{F}_{\lambda_k} = 
\begin{bmatrix}
\hat{Z}_{\lambda_k}^{(11)} & \cdots & \hat{Z}_{\lambda_k}^{(1s)} \\
\vdots & \ddots & \vdots \\
\hat{Z}_{\lambda_k}^{(s1)} & \cdots & \hat{Z}_{\lambda_k}^{(ss)}
\end{bmatrix}
\quad \mathrm{with} \quad
\hat{Z}_{\lambda_k}^{(ij)} = \begin{bmatrix} \beta_{i,j} + \mathbf{i} \gamma_{ij} \sqrt{\lambda_k} \\ & \beta_{i,j} - \mathbf{i} \gamma_{ij} \sqrt{\lambda_k} \end{bmatrix}.
$}
\end{equation*}
\noindent In words, we transformed $\tilde{Z}_{\lambda_k}$ by a similarity transformation into a block matrix where each block is diagonal and since these blocks are $2$-by-$2$, we can use the red-black symmetrical reordering of the unknowns represented by the permutation matrix $\tilde{\Pi} \in \mathbb{R}^{2s\times 2s}$ to obtain
\begin{equation}\label{eqn_secAnal_BlckDiagOf_tildeZlambdak}
\left(\tilde{F}_{\lambda_k} \tilde{\Pi} \right)^{-1} \tilde{Z}_{\lambda_k} \tilde{F}_{\lambda_k} \tilde{\Pi} = 
\begin{bmatrix} Z_{\lambda_k,+} \\ & Z_{\lambda_k,-} \end{bmatrix}
\end{equation}
\noindent with 
\begin{equation*}
Z_{\lambda_k,\pm} := \begin{bmatrix}
\beta_{11} \pm \mathbf{i} \gamma_{11} \sqrt{\lambda_k} & \cdots & \beta_{1s} \pm \mathbf{i} \gamma_{1s} \sqrt{\lambda_k} \\
\vdots & \ddots & \vdots \\
\beta_{s1} \pm \mathbf{i} \gamma_{s1} \sqrt{\lambda_k} & \cdots & \beta_{ss} \pm \mathbf{i} \gamma_{ss} \sqrt{\lambda_k}
\end{bmatrix} \in \mathbb{C}^{s\times s}.
\end{equation*}
\noindent The above calculations are similar to~\cite{outrata2022IRKILAS,outrata2023irksstage,outrata2025recent} and proves just as useful here. Having $\lambda_k^{(M,E)}$ (or their approximation), the eigenproperties of $\tilde{Z}_{\lambda_k}$ are governed by those of $Z_{\lambda_k,\pm}$ and eigenproblems with these small, $s$-by-$s$ matrices can be solved rapidly at virtually no cost. 

The last obstacle is the fact that~\eqref{eqn_secSpctrAnal_BlckDiagonalizationOfmtrxNNsZ_to2s} is not a similarity transformation and thus the eigenproperties of $Z_{\lambda_k,\pm}$ and $\tilde{Z}_{\lambda_k}$ don't relate to those of $\mtrxNNs{Z}$. However, we are interested in the eigenproperties of the \emph{preconditioned} system rather than those of the preconditioner and/or the system matrix. Repeating the calculation for the product $(\mtrxNNs{Z}^{\mathrm{prec}})^{-1}\mtrxNNs{Z}^{\mathrm{sys}}$ for some $\mtrxNNs{Z}^{\mathrm{prec}},\mtrxNNs{Z}^{\mathrm{sys}}$ of the form~\eqref{eqn_secAnal_Z_eq_bM_p_cE}, we obtain
\begin{equation}\label{eqn_secAnal_ZinvZ_BlockDiagonalized_w_2by2DiagBlcks}
\left( \mtrxNNs{Z}^{\mathrm{prec}}\right)^{-1} \mtrxNNs{Z}^{\mathrm{sys}} = 
\mtrxNNs{Q}
\begin{bmatrix}
\left( \mtrxN{Z}_{\lambda_1,+}^{\mathrm{prec}} \right)^{-1} \mtrxN{Z}_{\lambda_1,+}^{\mathrm{sys}} \\ 
& \ddots \\ && \left( \mtrxN{Z}_{\lambda_N,-}^{\mathrm{prec}} \right)^{-1}\mtrxN{Z}_{\lambda_N,-}^{\mathrm{sys}}
\end{bmatrix}
\mtrxNNs{Q}^{-1} ,
\end{equation}
\noindent with $\mtrxNNs{Q}:= \mtrxNNs{V} \Pi \mtrxNNs{F}$ and $\mtrxNNs{F}:= \mathrm{diag}\left( \tilde{F}_{\lambda_1} \tilde{\Pi},\tilde{F}_{\lambda_2} \tilde{\Pi}, \dotsc , \tilde{F}_{\lambda_N} \tilde{\Pi} \right)$. In words, the above calculation characterizes the eigenproperties of the preconditioned system $\mathcal{P}^{-1}\mathcal{T}$ in terms of $2N$ smaller eigenvalue problems of the size $s$-by-$s$. In particular, setting $\theta_k := \vec{i} h_t \sqrt{\lambda_k}$ and indexing with $\theta_k$ instead of $\lambda_k$ we have
\begin{equation*}
\resizebox{.995\textwidth}{!}{$
\mtrxN{T}_{\theta_k} = \begin{bmatrix}
1 - a_{11} \theta_k & - a_{12} \theta_k & \hdots & - a_{1s} \theta_k \\
- a_{21} \theta_k & 1 - a_{22} \theta_k & & \vdots \\
\vdots & & \ddots & \vdots \\
- a_{s1} \theta_k & \hdots & \hdots & 1 - a_{ss} \theta_k \\
\end{bmatrix}, \;
\mtrxN{P}_{\theta_k} = \begin{bmatrix}
1 - \alpha_{11} \theta_k & - \alpha_{12} \theta_k & \hdots & - \alpha_{1s} \theta_k \\
- \alpha_{21} \theta_k & 1 - \alpha_{22} \theta_k & & \vdots \\
\vdots & & \ddots & \vdots \\
- \alpha_{s1} \theta_k & \hdots & \hdots & 1 - \alpha_{ss} \theta_k \\
\end{bmatrix},
$}
\end{equation*}
\noindent the eigenproperties of the preconditioned system $\mathcal{P}^{-1}\mathcal{T}$ can be compacted into the following  parametric eigenvalue problem
\begin{equation}\label{eqn_secAnal_ParamEigPrblmFor_Xtheta}
\mtrxN{X}_{\pm \theta} \vec{s} = \xi \vec{s}
\quad \mathrm{with} \quad
\mtrxN{X}_{\theta} := \mtrxN{P}_{\theta}^{-1} \mtrxN{T}_{\theta}
\end{equation}
\noindent where $\theta \in \{ \theta_1, \dotsc , \theta_N \}$ and the eigenvectors $\vec{s}_k$ give the coefficients of the eigenvectors of $\mtrxNNs{P}^{-1}\mtrxNNs{T}$ in the basis given by the columns of $\mtrxNNs{Q}$; see~\cite[Section 3 and onward]{outrata2023irksstage}.

Assuming the pencil $\{M,E\}$ is real, the eigenvalues $\lambda_k$ come in complex conjugate pairs and so every $X_{\theta_k}$ can be paired with $X_{\bar{\theta}_{k}}$. Following a direct calculation (see~\cite[Proposition 3.7]{outrata2023irksstage}), we observe that the eigenpairs of $X_{\bar{\theta}_{k}}$ amount to the complex conjugation of those of $X_{\theta_{k}}$, i.e., we only need to consider one of these. Hence, we obtained the eigenproperties of the preconditioned system in terms of those of the pencil $\{M,E\}$ and a sequence of $s$-by-$s$ matrices $X_{\theta_k},\; (k=1,\dotsc , N)$.

We reached~\eqref{eqn_secAnal_ZinvZ_BlockDiagonalized_w_2by2DiagBlcks} assuming that $\lambda_k \neq 0$ for all $k$. If $\lambda_k = 0$ for some $k$, then $F_{\lambda_k}$ no longer exists: the matrix on the left-hand side of~\eqref{eqn_secAnal_ZinvZ_BlockDiagonalized_w_2by2DiagBlcks} becomes either a scaled identity matrix (if $\gamma = 0$) or a scaled Jordan cell (if $\gamma \neq 0$). In the first case, the calculation leading to~\eqref{eqn_secAnal_ZinvZ_BlockDiagonalized_w_2by2DiagBlcks} is immediately salvageable by taking $F_{\lambda_k} := I_2$ as the $2$-by-$2$ block is already in a diagonal form. In the second case we can apply the symmetrical permutation $\tilde{\Pi}$ to $\tilde{Z}_{\lambda_k}$ as in~\eqref{eqn_secAnal_BlckDiagOf_tildeZlambdak} and obtain
\begin{equation*}
\tilde{\Pi}^T \tilde{Z}_{\lambda_k} \tilde{\Pi} = \tilde{\Pi}^T \left( B \otimes I_{2} + \Gamma \otimes J_2 \right) \tilde{\Pi} = I_2 \otimes B + J_2 \otimes \Gamma ,
\end{equation*}
\noindent with
\begin{equation*}
B := 
\begin{bmatrix}
\beta_{11} & \cdots & \beta_{1s} \\
\vdots & \ddots & \vdots \\
\beta_{s1} & \cdots & \beta_{ss} 
\end{bmatrix},  \quad
\Gamma :=
\begin{bmatrix}
\gamma_{11} & \cdots & \gamma_{1s} \\
\vdots & \ddots & \vdots \\
\gamma_{s1} & \cdots & \gamma_{ss}
\end{bmatrix}
\quad \mathrm{and} \quad
J_2 = \begin{bmatrix} 0 & 1 \\ 0 & 0 \end{bmatrix}.
\end{equation*}
\noindent Recalling that in our specific setting we have $B=I_s$ and $\Gamma=A$ for $\mtrxNNs{T}$ (or $\Gamma=\tilde{A}$ for $\mtrxNNs{P}$), we see that the if $\lambda_k=0$ then we replace the two consecutive diagonal blocks $\left( \mtrxN{Z}_{\lambda_k,+}^{\mathrm{prec}} \right)^{-1} \mtrxN{Z}_{\lambda_k,+}^{\mathrm{sys}}$ and $\left( \mtrxN{Z}_{\lambda_k,-}^{\mathrm{prec}} \right)^{-1}\mtrxN{Z}_{\lambda_k,-}^{\mathrm{sys}}$ in~\eqref{eqn_secAnal_ZinvZ_BlockDiagonalized_w_2by2DiagBlcks} by
\begin{equation}\label{eqn_secAnal_CalcOfDiagonalizabilityForLambdaeq0_Atild_m_A}
\begin{bmatrix} I & -h_t\tilde{A} \\ & I \end{bmatrix}^{-1} \begin{bmatrix} I & -h_t A \\ & I \end{bmatrix} = 
\begin{bmatrix} I & h_t (\tilde{A} - A) \\ & I \end{bmatrix} \in \mathbb{R}^{2s\times 2s}.
\end{equation}
\noindent This shows that such case results in $2s$ eigenvalues of the preconditioned system equal to $1$ but the number of eigenvectors amounts only to $2s - \mathrm{rank}(\tilde{A} - A)$, i.e., each $\lambda_k=0$ results in $\mathrm{rank}(\tilde{A} - A)$ of Jordan chains of length $2$. Notably, the eigenvalues are included in the limit of $\theta_k \rightarrow 0$ as the eigenvalues are continuous functions of matrix entries.

Assuming $\lambda_k\neq 0$ for all $k$, if $s>3$ a direct calculation of the spectrum and eigenvectors of $\mtrxN{X}_{\theta_k}$ becomes unwieldy or downright intractable and so let us assume that $\mtrxN{X}_{\theta_k}$ diagonalizes under some similarity transformation, denoted by $\mtrxN{S}_{\theta_k}$, 
\begin{equation}\label{eqn_secSpctrAnalKroneckr_XkronTheta_eq_SkronThetaXikronThetaSkronThetaInv}
\mtrxN{X}_{\theta_k} = \mtrxN{S}_{\theta_k} \begin{bmatrix}\xi_{\theta_k}^{(1)} \\ & \ddots \\ && \xi_{\theta_k}^{(s)}\end{bmatrix} \mtrxN{S}_{\theta_k}^{-1}.
\end{equation}
\noindent Due to~\eqref{eqn_secAnal_ParamEigPrblmFor_Xtheta}, the spectrum of the preconditioned system $\mtrxNNs{P}^{-1}\mtrxNNs{T}$ corresponds to $\{ \xi_{\theta_k}^{(i)}\}_{i,k}$. As for the condition number of the eigenbasis $\mtrxNNs{S}$ of $\mtrxNNs{P}^{-1}\mtrxNNs{T}$ we can bound it by 
\begin{equation}\label{eqn_secEigbase_BoundForCondNmbEigBasisSkron}
\kappa \left( \mtrxNNs{S} \right) \leq 
\kappa ( \mtrxNNs{V} ) \cdot 
\left( \max\limits_{k=1,\dotsc ,N} \| F_{\theta_k} \| \right) \left( \max\limits_{k=1,\dotsc ,N} \| F_{\theta_k}^{-1} \| \right)
\cdot 
\left( \max\limits_{k=1,\dotsc ,N} \| S_{\theta_k} \| \right) \left( \max\limits_{k=1,\dotsc ,N} \| S_{\theta_k}^{-1} \| \right),
\end{equation}
\noindent and thanks to the direct computation
\begin{equation*}
\| F_{\theta_k} \|_2 = \sqrt{2} \max \{ 1, 1/|\sqrt{\lambda_k}| \}
\quad \mathrm{and} \quad
\| F_{\theta_k}^{-1} \|_2 = 1/\sqrt{2} \max \{ 1, |\sqrt{\lambda_k}| \}
 \end{equation*}
\noindent we observe that the product of the second and third terms in~\eqref{eqn_secEigbase_BoundForCondNmbEigBasisSkron} scales as 
\begin{equation}\label{eqn_secSpctrAnalKroneckr_normF_AndInv_scaling}
\max \{ 1, \sqrt{|\lambda_{\mathrm{max}}|} , 1/\sqrt{|\lambda_{\mathrm{min}}|}, \sqrt{|\lambda_{\mathrm{max}}|}/\sqrt{|\lambda_{\mathrm{min}}|} \},
\end{equation}
\noindent where $\lambda_{\mathrm{max, min}}$ are the extremal generalized eigenvalues of $\{M,E\}$ in magnitude.

\begin{claim}
The above calculations provide the required mathematical background for understanding and estimating the observed numerical results in~\cite{rani2025efficient} for the benchmark preconditioners $\mtrxNNs{P}_{\kappa, LD, TAI}$.
\end{claim}

The two main observations for these preconditioners are about the conditioning of the eigenbasis of the preconditioned system and the average number of GMRES iterations, see~\cite[Table 4.1 and Figure 4.8]{rani2025efficient}. First,~\cite[Table 4.1]{rani2025efficient} complements the spectral graphs of the preconditioners and links these to GMRES convergence bounds and make several claims about the data without providing explanation -- we deliver on the promise of providing the explanation below. Second, showing that the number of GMRES iterations is expected to stay constant, roughly around the reported number of iterations in~\cite[Figure 4.8]{rani2025efficient}, together with~\cite[Section 4.3]{rani2025efficient} firmly grounds the observed runtimes under mesh refinement. For both of these tasks, we consider here only the preconditioned system $\mtrxNNs{P}_{TAI}^{-1}\mtrxNNs{T}$; the other two results are qualitatively analogous (while quantitatively different) and the interested reader is welcomed to explore these using the provided code at \href{https://github.com/MichalOutrata/PrecondAnalysisForIRK/}{https://github.com/MichalOutrata/PrecondAnalysisForIRK/}.

\section{The eigenbasis of the preconditioned system}\label{sec_CondOfEigBasis}

First, the results in~\cite[Table 4.1]{rani2025efficient} correspond to $\mathscr{L}=c^2\Delta$ in~\eqref{eqn_secIntro_PDE_utt_m_L_eq_f} with \emph{Neumann boundary conditions}, i.e., the spatial problem is \emph{singular}: constant functions lie in the kernel. The stiffness matrix $\mtrxN{E}$ inherits this property, having one constant singular mode. As we never solve any systems with $\mtrxN{E}$ this doesn't pose any problems but the above analysis shows a single zero eigenvalue makes the system not diagonalizable: the matrix $\tilde{A}-A$ is non-zero even though it is not always full rank\footnote{In our experience, for the $TAI, LD$ approaches, the matrix $\tilde{A}-A$ has usually a one-dimensional kernel while for the $\kappa$ approach the matrix $\tilde{A}-A$ is usually regular.}, see~\eqref{eqn_secAnal_CalcOfDiagonalizabilityForLambdaeq0_Atild_m_A}. Numerically, we therefore obtain $1/\sqrt{\lambda_{\mathrm{min}}} \approx \varepsilon_{\mathrm{mach}}^{-1/2}$ and get the conditioning of the eigenbasis of the preconditioned system of approximately $10^8$, as reported in~\cite[Table 4.1]{rani2025efficient} -- in spite of the matrix being nondiagonalizable! This also explains the independence with respect to the discretization parameters $s, h_x$ and $h_t$. Moreover, this insight immediately disqualifies the use of the GMRES bound~\eqref{eqn_secSolv_GMRESPolyInEigenvals} (and its descendants), although computing $\kappa(\mtrxNNs{S})$ numerically wouldn't suggest so.

Second, having $\lambda_1=0$, the correct analogue of~\eqref{eqn_secAnal_ZinvZ_BlockDiagonalized_w_2by2DiagBlcks}  becomes
\begin{equation}\label{eqn_secCondOfEigBasis_PinvT_eq_BlckDiag}
\mtrxNNs{P}^{-1} \mtrxNNs{T} = 
\mtrxNNs{Q}
\begin{bmatrix}
\mtrxN{I} & h_t(\tilde{A}-A) \\
& \mtrxN{I} \\ 
&& \mtrxN{X}_{\theta_2} \\
&&& \ddots \\
&&&& \mtrxN{X}_{-\theta_N}
\end{bmatrix}
\mtrxNNs{Q}^{-1} ,
\quad \mathrm{with} \quad
\mtrxN{X}_{\theta} := \mtrxN{P}_{\theta}^{-1} \mtrxN{T}_{\theta}.
\end{equation}
\noindent This shows that the nondiagonalizability of the preconditioned system due to $\lambda_1=0$ is inherently ``small'' and could be removed. Naturally, a particular choice of $A$ and $\tilde{A}$ can result in a nondiagonalizable $\mtrxN{X}_{\theta_k}$ for some $k$ but for the practical choices based on Gauss-Legendre, Radau, Lobatto IRK methods this hasn't been generally the case. Trying to estimate the GMRES behavior, a heuristic idea is to simply separate the invariant subspaces of the preconditioned system induced by $\lambda_1$ and $\lambda_2,\dotsc ,\lambda_N$ and use the bound~\eqref{eqn_secSolv_GMRESPolyInEigenvals} separately, as if for the system matrix induced by $\lambda_2,\dotsc ,\lambda_N$. While this is mathematically incorrect for obtaining a bound, it can sometimes provide an intuition for the GMRES behavior; the bound on the conditioning of the eigenbases when excluding the subspace induced by $\lambda_1=0$ are the numbers in the parenthesis in the first three columns in~\cite[Table 4.1]{rani2025efficient} and these scale proportionally to the bound~\eqref{eqn_secEigbase_BoundForCondNmbEigBasisSkron}.

Last, we want to expand on the comment accompanying~\cite[Table 4.1]{rani2025efficient} ``\emph{\dots to obtain these results accurately, requires a careful reformulation of the calculation \dots}''. We ran numerous numerical experiments using the black-box routines\footnote{We tried both \texttt{np.linalg.eig()},\texttt{np.linalg.cond()} in python and \texttt{eig} and \texttt{cond} in MATLAB, all relying on LAPACK's \texttt{\_geev()}.} for computing the spectral factorization of $\mtrxNNs{P}^{-1}\mtrxNNs{T}$ and then for the condition number of the resulting eigenvector matrix. The results usually \emph{vastly} overestimate the bound of the condition number $\kappa(\mtrxNNs{S})$ in~\eqref{eqn_secEigbase_BoundForCondNmbEigBasisSkron}. First, we acknowledge that this is not too surprising for the problem at hand in~\cite[Table 4.1]{rani2025efficient} as we have established it is not diagonalizable. However, even when considering $\mathcal{L}$ that is non-singular, this huge overestimation remained true, often by \emph{many} orders of magnitude, see the mentioned codebase for illustrations. Only after we permute the preconditioned system into a block-diagonal one, like in~\eqref{eqn_secCondOfEigBasis_PinvT_eq_BlckDiag}, the correct values are retained. After consulting our experience with colleagues, our intuition is that the reason is the parametric nature of the eigenproperties of the preconditioned system -- clearly visible in the block-diagonal form but obscured in the naive form: in all of our test problems, the eigenvalues $\lambda_1,\dotsc ,\lambda_N$ populate some interval (and they populate it fairly densely in certain parts) and the parametrized eigenvalues problem then corresponds to (fairly dense) sampling over that interval. While eigenvalues are continuous with respect to the matrix entries, eigenvectors are not and this set-up poses issues for accurate resolution of the eigenvectors for the standard black-box solvers; they need the explicit structural information to be able to resolve the eigenvectors accurately.

\section{The preconditioned GMRES behavior}\label{sec_precGMRESbeh}
We want to put the eigenproperties of $\mtrxNNs{P}^{-1}\mtrxNNs{T}$ to a good use and obtain a cheap convergence estimate of~\eqref{eqn_secSolv_GMRESAsymptConvEst}, i.e., we need to \emph{choose} and \emph{justify} $\sigma^{non-discr}$ for which we can compute $\rho_{\mathrm{est}}$ and $\kappa_{\mathrm{est}}$. Both of these points were treated in~\cite[Section 3.1]{outrata2023irksstage} for a parabolic problem instead of~\eqref{eqn_secIntro_PDE_utt_m_L_eq_f} but the same idea carries over.

\begin{figure}[t]
\centering
\resizebox{.99\textwidth}{!}{
	\includegraphics{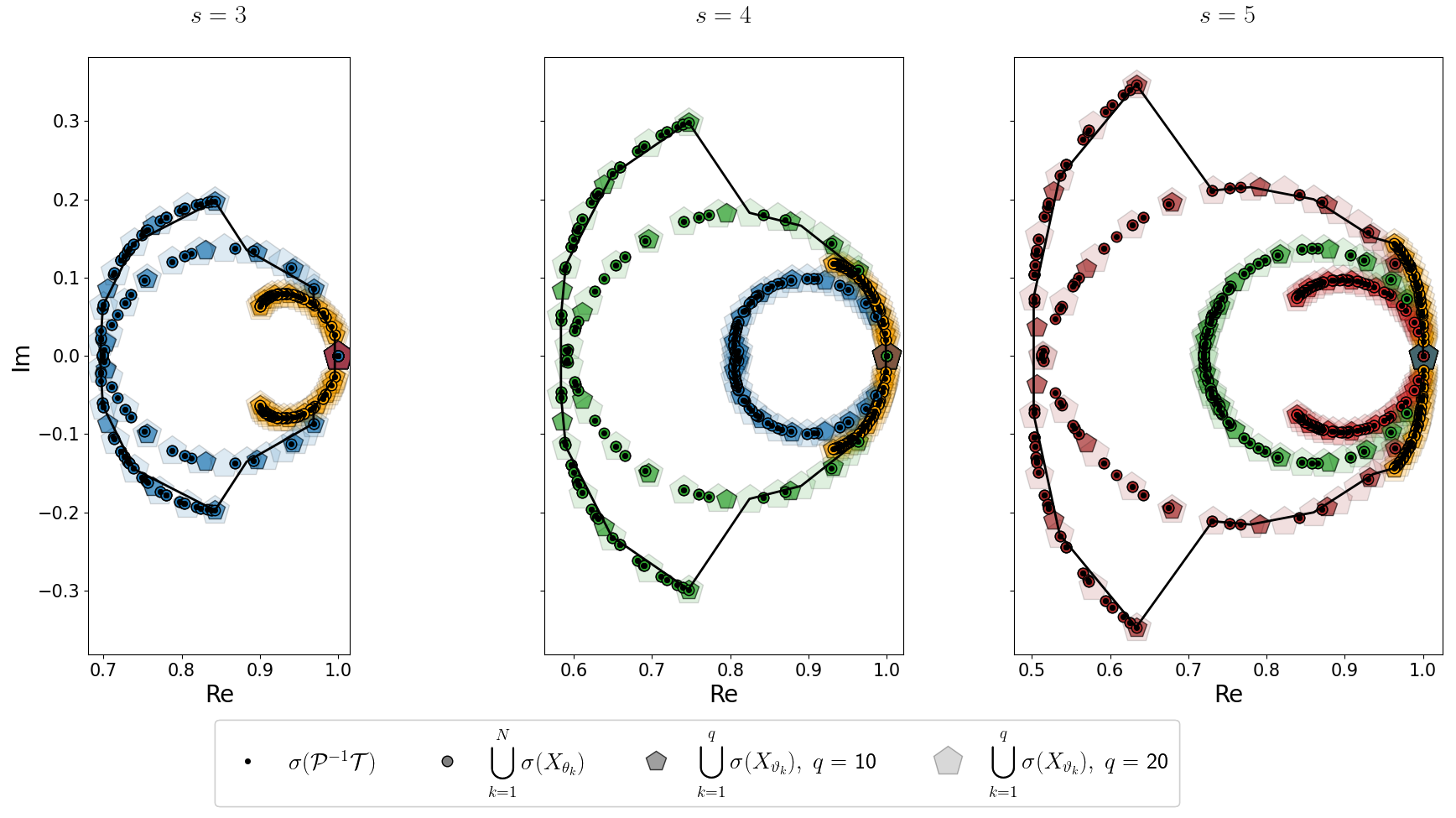}
}
\caption{For the TAI preconditioner, we show the spectrum of the assembled preconditioned system $\sigma(\mtrxNNs{P}^{-1}\mtrxNNs{T})$, its reformulation $\cup_{k=1}^{N} \sigma(X_{\theta_k})$ and its approximations using sampling with $\vartheta_1,\dotsc , \vartheta_q \in \mathcal{C}_{\theta}$ for $q=10,20$. Different branches of the underlying algebraic curve $\Gamma$ are highlighted in different colors. We used the set-up of~\cite[Figures 4.2 and 4.3]{rani2025efficient}, i.e., $\mathcal{L} = \Delta$ on the unit square with Neumann BC with $N=81$ (for larger $N$ the different markers and colors are hard to discern, the process itself works for much larger $N$, see the mentioned github repository) and Gauss-Legendre IRK.}\label{fig_secAnal_SpectraAsCurves}
\end{figure}

The first important observation is that the set $\{ \theta_k \}_{k=1,\dotsc ,N}$ is a simple, smooth transformation of the generalized spectrum of the pencil $\{\mtrxN{M},\mtrxN{E}\}$ and that the eigenvalue problem for the preconditioned system matrix $\mtrxNNs{P}^{-1}\mtrxNNs{T}$ can be written as a sampling of a particular \emph{parametrized eigenvalue problem} precisely at the points $\{ \theta_k \}_{k}$. We will assume that as the spatial discretization refines, the set $\{ \theta_k \}_{k}$ starts populating densely a curve $\mathcal{C}_{\theta} \subset \mathbb{C}$, e.g., $\mathcal{C}_{\theta} = \vec{i}(0,+\infty)$. Denoting the eigenvalue curves of the parametrized eigenvalues problems over $\mathcal{C}_{\theta}$ as $\Gamma_1, \dotsc , \Gamma_s$, it seems natural to take $\sigma^{non-discr}= \Gamma := \Gamma_1 \, \cup \, \dotsc  \, \cup \, \Gamma_s$ -- it is \emph{the} asymptotically correct choice with respect to mesh refinement in space. However, in practice, the eigenvalue curves $\Gamma_1, \dotsc , \Gamma_s$ are not analytically tractable. We instead sample\footnote{In case $\mathcal{C}_{\theta}$ is unbounded like in~\cite{rani2025efficient}, we need to estimate the endpoints based on the spectrum of $\mtrxNNs{L}$.} $\mathcal{C}_{\theta}$ with few \emph{artificial} values $\vartheta_1, \dotsc , \vartheta_q$, obtaining $q$ points on each $\Gamma_1\dotsc , \Gamma_s$ and interpolating these we get $\sigma^{non-discr}:= \Gamma_{\vartheta} \approx \Gamma$, i.e., we interpolate the ``asymptotic spectrum of $\mtrxNNs{P}^{-1}\mtrxNNs{T}$'' based on the eigenvalues of $\{ X_{\vartheta_k} \}_{k=1,\dotsc ,q}$. We illustrate this in Figure~\ref{fig_secAnal_SpectraAsCurves} for two particular choices of preconditioners for a fairly small problem.

Notably, in~\cite{rani2025efficient}, $\Gamma$ stays uniformly bounded, even if $\mathcal{C}_{\theta}$ is unbounded as
\begin{equation*}
\lim\limits_{\theta \rightarrow 0} X_{\theta} = I
\quad \mathrm{and} \quad
\lim\limits_{ | \theta | \rightarrow +\infty } X_{\theta} \rightarrow \tilde{A}^{-1}A.
\end{equation*}

The second important observation is that for compact and connected sets the min-max problem in~\eqref{eqn_secSolv_GMRESPolyOverSigmaNonDiscr} can be asymptotically evaluated using \emph{potential theory} and \emph{Schwarz-Christoffel maps}, see~\cite[Sections 2 and 3]{driscoll1998potential}. This process uses the so-called \emph{exterior map of $\sigma^{non-discr}$}, see~\cite{driscoll2002schwarz}; numerically this can be done using the excellent \texttt{sctoolbox} in MATLAB~\cite{driscoll1996algorithm}, assuming $\sigma^{non-discr}$ is a polygonal domain or curve\footnote{If $\sigma^{non-discr}$ is a union of polygonal domains (or possibly line segments) that are symmetric with respect to the real axis, the authors in~\cite{embree1999green} reformulate the problem of calculating the corresponding $\rho_{\mathrm{est}}$ so that the \texttt{sctoolbox} can be used again, requiring, however, more preliminary work for the reformulation.}. Here, as in~\cite[Sections 3 and 4]{outrata2023irksstage}, we take an ``\emph{envelope}'' that tightly encloses $\sigma^{non-discr} \equiv \Gamma_{\vartheta}$ but is simply connected and without isolated points (showed on Figure~\ref{fig_secAnal_SpectraAsCurves}) as the black line). 

As for $\kappa_{\mathrm{est}}$, the artificially sampled points $\vartheta_1,\dotsc ,\vartheta_q$ can be used to estimate
$\kappa( \mtrxNNs{S})$: instead of calculating only the eigenvalues of $X_{\vartheta_k}$, we calculate their eigenbasis too and use~\eqref{eqn_secEigbase_BoundForCondNmbEigBasisSkron} to obtain an estimate on $\kappa( \mtrxNNs{S})$ (assuming these are diagonalizable). We keep the same setting as in Figure~\ref{fig_secAnal_SpectraAsCurves} and we see that already for a very small $q$, the estimates capture the GMRES convergence behavior accurately.

\begin{figure}[t]
\centering
\resizebox{.9\textwidth}{!}{
	\includegraphics{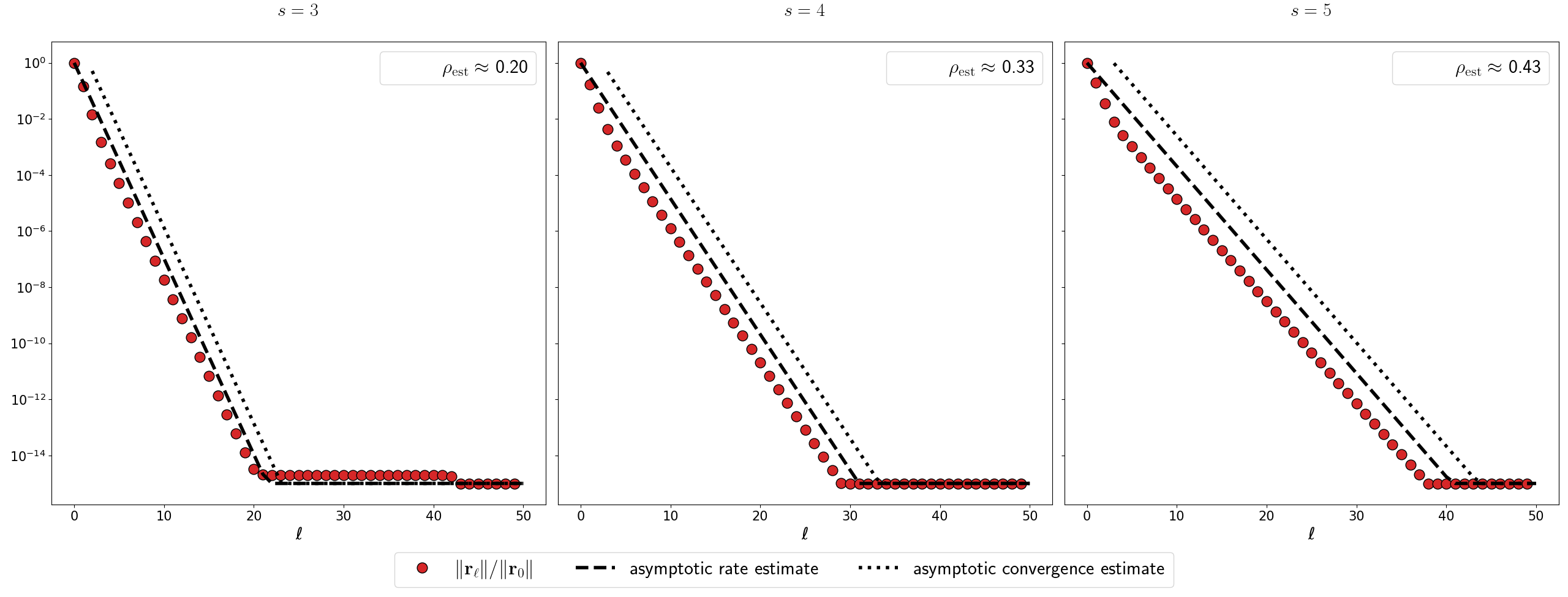}
}
\caption{We show the GMRES convergence behavior and its estimates for the TAI preconditioner and the set-up of~\cite[Figures 4.2 and 4.3]{rani2025efficient}, i.e., $\mathcal{L} = \Delta$ on the unit square with Neumann BC and Gauss-Legendre IRK, when neglecting the singular mode of the spatial operator for the estimates (otherwise the $\kappa_{\mathrm{est}}$ will correspond to roughly the values in~\cite[Table 4.1]{rani2025efficient}) as explained in Section~\ref{sec_CondOfEigBasis}). We took $N=16641$ and obtained the estimates with $q=20$.}\label{fig_secAnal_GMRES_EstAndSpctr_PinT}
\end{figure}

\bibliographystyle{plainnat}
\bibliography{biblio}

@InProceedings{Greenbaum1994a,
  author    = {Greenbaum, A. and Strako{\v{s}}, Z.},
  booktitle = {{R}ecent {A}dvances in {I}terative {M}ethods},
  title     = {Matrices that generate the same {K}rylov residual spaces},
  year      = {1994},
  editor    = {Golub, G. H. and Greenbaum, A. and Luskin, M.},
  pages     = {95--118},
  publisher = {Springer},
  series    = {{IMA} {V}olumes in {M}athematics and its {A}pplications},
  volume    = {60},
}

@Article{Greenbaum1996,
  author  = {Greenbaum, A. and Pt{\'a}k, V. and Strako{\v{s}}, Z.},
  journal = {{SIAM} {J}ournal on {M}atrix {A}nalysis and {A}pplications},
  title   = {Any nonincreasing convergence curve is possible for {GMRES}},
  year    = {1996},
  number  = {3},
  pages   = {465--469},
  volume  = {17},
}

@Article{Arioli1998,
  author  = {Arioli, M. and Pt{\'a}k, V. and Strako{\v{s}}, Z.},
  journal = {{BIT} {N}umerical {M}athematics},
  title   = {Krylov sequences of maximal length and convergence of {GMRES}},
  year    = {1998},
  number  = {4},
  pages   = {636--643},
  volume  = {38},
}

@Book{liesen2013krylov,
  author    = {Liesen, J. and Strako{\v{s}}, Z.},
  publisher = {Oxford {U}niversity {P}ress},
  title     = {Krylov {S}ubspace {M}ethods: {P}rinciples and {A}nalysis},
  year      = {2012},
  doi       = {10.1093/acprof:oso/9780199655410.001.0001},
}

@InProceedings{Neytcheva2020,
  author       = {Neytcheva, M. and Axelsson, O.},
  booktitle    = {{P}roceedings of the {C}onference {A}lgoritmy 2020},
  title        = {Numerical {S}olution {M}ethods for {I}mplicit {R}unge-{K}utta {M}ethods of {A}rbitrarily {H}igh {O}rder},
  year         = {2020},
  editor       = {Frolkovi{\v{c}}, P. and Mikula, K. and {\v{S}}ev{\v{c}}ovi{\v{c}}, D.},
  organization = {{S}lovak {U}niversity of {T}echnology in {B}ratislava},
  publisher    = {Vydavate{\'{l}}stvo {SPEKTRUM}},
}

@Article{Howle2021new,
  author  = {Rana, M. M. and Howle, V. E. and Long, K. and Meek, A. and Milestone, W.},
  journal = {{SIAM} {J}ournal on {S}cientific {C}omputing},
  title   = {A new block preconditioner for implicit {R}unge-{K}utta methods for parabolic {PDE} problems},
  year    = {2021},
  number  = {5},
  pages   = {S475--S495},
  volume  = {43},
}

@Article{southworth2022fastI,
  author    = {Southworth, B. S. and Krzysik, O. and Pazner, W. and De Sterck, H.},
  journal   = {{SIAM} {J}ournal on {S}cientific {C}omputing},
  title     = {Fast solution of fully implicit {R}unge--{K}utta and discontinuous {G}alerkin in time for numerical {PDE}s, {P}art {I}: {T}he linear setting},
  year      = {2022},
  number    = {1},
  pages     = {416--443},
  volume    = {44},
  publisher = {{SIAM}},
}

@Article{southworth2022fastII,
  author    = {Southworth, B. S. and Krzysik, O. and Pazner, W.},
  journal   = {{SIAM} {J}ournal on {S}cientific {C}omputing},
  title     = {Fast solution of fully implicit {R}unge--{K}utta and discontinuous {G}alerkin in time for numerical {PDE}s, {P}art {II}: nonlinearities and {DAE}s},
  year      = {2022},
  number    = {2},
  pages     = {636--663},
  volume    = {44},
  publisher = {{SIAM}},
}

@Book{bai2000templates,
  author    = {Bai, Z. and Demmel, J. and Dongarra, J. and Ruhe, A. and van der Vorst, H.},
  publisher = {{SIAM}, {P}hiladelphia},
  title     = {Templates for the {S}olution of {A}lgebraic {E}igenvalue {P}roblems: {A} {P}ractical {G}uide},
  year      = {2000},
}

@Article{staff2006preconditioning,
  author    = {Staff, G. A. and Mardal, K.-A. and Nilssen, T. K.},
  journal   = {Modeling, {I}dentification and {C}ontrol},
  title     = {{Preconditioning of fully implicit {R}unge-{K}utta schemes for parabolic {PDE}s}},
  year      = {2006},
  number    = {2},
  pages     = {109--123},
  volume    = {27},
  publisher = {Norwegian {S}ociety of {A}utomatic {C}ontrol},
}

@Article{outrata2022IRKILAS,
  author  = {Gander, M. J. and Outrata, M.},
  journal = {{L}inear {A}lgebra and its {A}pplications},
  title   = {Spectral analysis of implicit 2-stage block {R}unge-{K}utta preconditioners},
  year    = {2023},
  issn    = {0024-3795},
  doi     = {10.1016/j.laa.2023.07.008},
}

@Misc{howle2022efficient,
  author        = {Clines, M. R. and Howle, V. E. and Long, K. R.},
  title         = {Efficient order-optimal preconditioners for implicit {R}unge-{K}utta and {R}unge-{K}utta-{N}ystr{\"{o}}m methods applicable to a large class of parabolic and hyperbolic {PDE}s},
  year          = {2022},
  archiveprefix = {arXiv},
  eprint        = {2206.08991},
}

@Article{driscoll1998potential,
  author    = {Driscoll, T. A. and Toh, K.-C. and Trefethen, L. N.},
  journal   = {{SIAM} {R}eview},
  title     = {From potential theory to matrix iterations in six steps},
  year      = {1998},
  number    = {3},
  pages     = {547--578},
  volume    = {40},
  publisher = {{SIAM}},
}

@TechReport{driscoll1996algorithm,
  author    = {Driscoll, T. A.},
  title     = {A {MATLAB} toolbox for {S}chwarz-{C}hristoffel mapping},
  year      = {1996},
  number    = {2},
  journal   = {{ACM} {T}ransactions on {M}athematical {S}oftware},
  pages     = {168--186},
  publisher = {{ACM} {N}ew {Y}ork, {NY}, {USA}},
  volume    = {22},
}

@Book{driscoll2002schwarz,
  author    = {Driscoll, T. A. and Trefethen, L. N.},
  publisher = {Cambridge {U}niversity {P}ress, {C}ambridge},
  title     = {Schwarz-{C}hristoffel mapping},
  year      = {2002},
  edition   = {{F}irst},
}

@Article{embree1999green,
  author    = {Embree, M. and Trefethen, L. N.},
  journal   = {{SIAM} {R}eview},
  title     = {Green's functions for multiply connected domains via conformal mapping},
  year      = {1999},
  number    = {4},
  pages     = {745--761},
  volume    = {41},
  publisher = {{SIAM}},
}

@Article{outrata2023irksstage,
  author  = {Gander, M. J. and Outrata, M.},
  journal = {{SIAM} {J}ournal on {S}cientific {C}omputing, in press},
  title   = {Spectral analysis of implicit $s$-stage block {R}unge-{K}utta preconditioners},
  year    = {2024},
  number  = {3},
  pages   = {A2047-A2072},
  volume  = {46},
  doi     = {10.1137/23M1604266},
}

@Article{dravins2024stageparallel_2,
  author  = {Munch, P. and Dravins, I. and Kronbichler, M. and Neytcheva, M.},
  journal = {{SIAM} {J}ournal on {S}cientific {C}omputing},
  title   = {Stage-parallel fully implicit {R}unge–{K}utta implementations with optimal multilevel preconditioners at the scaling limit},
  year    = {2024},
  number  = {2},
  pages   = {S71-S96},
  volume  = {46},
  doi     = {10.1137/22M1503270},
}

@Article{dravins2024spectral,
  author  = {Dravins, I. and Serra-Capizzano, S. and Neytcheva, M.},
  journal = {{SIAM} {J}ournal on {M}atrix {A}nalysis and {A}pplications},
  title   = {Spectral analysis of preconditioned matrices arising from stage-parallel implicit {R}unge–{K}utta methods of arbitrarily high order},
  year    = {2024},
  number  = {2},
  pages   = {1007-1034},
  volume  = {45},
  doi     = {10.1137/23M1552498},
}

@Article{dravins2024stageparallel_1,
  author  = {Axelsson, O. and Dravins, I. and Neytcheva, M.},
  journal = {{N}umerical {L}inear {A}lgebra with {A}pplications},
  title   = {Stage-parallel preconditioners for implicit {R}unge–{K}utta methods of arbitrarily high order, linear problems},
  year    = {2024},
  number  = {1},
  pages   = {e2532},
  volume  = {31},
  doi     = {10.1002/nla.2532},
}

@Article{butcher1976implementation,
  author    = {Butcher, J. C.},
  journal   = {{BIT} {N}umerical {M}athematics},
  title     = {On the implementation of implicit {R}unge-{K}utta methods},
  year      = {1976},
  number    = {3},
  pages     = {237--240},
  volume    = {16},
  doi       = {10.1007/BF01932265},
  publisher = {Springer},
}

@Article{pazner2017stage,
  author  = {Pazner, W. and Persson, P.-O.},
  journal = {{J}ournal of {C}omputational {P}hysics},
  title   = {Stage-parallel fully implicit {Runge–Kutta} solvers for discontinuous {G}alerkin fluid simulations},
  year    = {2017},
  pages   = {700-717},
  volume  = {335},
  doi     = {10.1016/j.jcp.2017.01.050},
}

@Article{pearson2024pint,
  author  = {Leveque, S. and Bergamaschi, L. and Mart\'{\i}nez, \'{A}. and Pearson, J. W.},
  journal = {{SIAM} {J}ournal on {M}atrix {A}nalysis and {A}pplications},
  title   = {Parallel-in-Time Solver for the All-at-Once {Runge–Kutta} Discretization},
  year    = {2024},
  number  = {4},
  pages   = {1902-1928},
  volume  = {45},
  doi     = {10.1137/23M1567862},
}

@Article{farrell2021irksome,
  author    = {Farrell, P. E. and Kirby, R. C. and Marchena-Men{\'{e}}ndez, J.},
  journal   = {{ACM} {T}ransactions on {M}athematical {S}oftware},
  title     = {Irksome: Automating {Runge–Kutta} Time-stepping for Finite Element Methods},
  year      = {2021},
  number    = {4},
  volume    = {47},
  address   = {{N}ew {Y}ork, {NY, USA}},
  doi       = {10.1145/3466168},
  numpages  = {26},
  publisher = {Association for {C}omputing {M}achinery},
}

@Article{rani2025efficient,
  author  = {Rani, A. and Ghysels, P. and Howle, V. and Long, K. and Outrata, M.},
  journal = {{SIAM} {J}ournal on {S}cientific {C}omputing},
  title   = {Efficient Solution of Fully Implicit {R}unge–{K}utta Methods for Linear Wave Equations},
  year    = {2025},
  number  = {0},
  pages   = {S183-S206},
  volume  = {0},
  doi     = {10.1137/24M1677484},
}

@Misc{dravins2024performance,
  author        = {Dravins, I. and Koch, M. and Griehl, V. and Kormann, K.},
  title         = {Performance evaluation of mixed-precision {R}unge-{K}utta methods for the solution of partial differential equations},
  year          = {2024},
  archiveprefix = {arXiv},
  eprint        = {2412.16638},
}

@Article{outrata2025recent,
  author  = {Outrata, M.},
  journal = {Proceedings in {A}pplied {M}athematics and {M}echanics},
  title   = {On Recent Advances of Spectral Analysis for Systems Arising From Fully-Implicit {RK} Methods},
  year    = {2026},
  number  = {1},
  pages   = {e70082},
  volume  = {26},
  doi     = {10.1002/pamm.70082},
}

\end{document}